\documentclass[11pt]{article}
\usepackage{graphicx}
\usepackage{epsfig}
\usepackage{latexsym}
\newsavebox{\toy}
\savebox{\toy}{\framebox[0.65em]{\rule{0cm}{1ex}}}
\newcommand{\QED}{\usebox{\toy}\end{demo}}
\newenvironment{property}%
{\begin{list}{}{\setlength{\rightmargin}{0pt}%
\setlength{\itemsep}{0pt}}}{\end{list}}
\newlength{\templength}
\newcommand{\bp}{\setlength{\templength}{\labelwidth}%
\setlength{\labelwidth}{2em}\begin{property}}
\newcommand{\ep}{\end{property}\setlength{\labelwidth}{\templength}}
\newtheorem{theorem}{Theorem}[subsection]
\newtheorem{lemma}[theorem]{Lemma}
\newtheorem{proposition}[theorem]{Proposition}
\newtheorem{corollary}[theorem]{Corollary}
\newtheorem{assumption}{Assumption}
\newtheorem{definition}{Definition}[subsection]
\newtheorem{remark}{Remark}[subsection]
\newtheorem{exercise}{Exercise}[subsection]

\newcommand{\Thm}[1]{Theorem \ref{Thm.#1}}
\newcommand{\Lem}[1]{Lemma \ref{Lem.#1}}

\newcommand{\Prop}[1]{Proposition \ref{Prop.#1}}

\newcommand{\Theorem}[1]{\begin{theorem}\label{Thm.#1}}
\newcommand{\Lemma}[1]{\begin{lemma}\label{Lem.#1}}
\newcommand{\Proposition}[1]{\begin{proposition}\label{Prop.#1}}
\newcommand{\Corollary}[1]{\begin{corollary}\label{Cor.#1}}
\newcommand{\Assumption}[1]{\begin{assumption}\label{Ass.#1}\rm}
\newcommand{\Definition}[1]{\begin{definition}\label{Def.#1}\rm}
\newcommand{\Remark}[1]{\begin{remark}\label{Rem.#1}\rm }
\newcommand{\Exercise}[1]{\begin{exercise}\label{Exe.#1}\rm }

\newcommand{\bd}{\begin{displaymath}}
\newcommand{\ed}{\end{displaymath}}
\newcommand{\bdn}{\begin{equation}}
\newcommand{\bdnl}{\begin{equation}\label}
\newcommand{\edn}{\end{equation}}
\newcommand{\barray}{\begin{array}}
\newcommand{\earray}{\end{array}}
\newcommand{\bds}{\begin{description}}
\newcommand{\eds}{\end{description}}
\newcommand{\bitemize}{\begin{itemize}}
\newcommand{\eitemize}{\end{itemize}}
\newcommand{\benumerate}{\begin{enumerate}}
\newcommand{\eenumerate}{\end{enumerate}}
\newcommand{\btabbing}{\begin{tabbing}}
\newcommand{\etabbing}{\end{tabbing}}
\newcommand{\bcenter}{\begin{center}}
\newcommand{\ecenter}{\end{center}}
\newcommand{\bflushright}{\begin{flushright}}
\newcommand{\bflushleft}{\begin{flushleft}}
\newcommand{\eflushright}{\end{flushright}}
\newcommand{\eflushleft}{\end{flushleft}}
\newcommand{\bdnn }{\begin{eqnarray*}}
\newcommand{\ednn }{\end{eqnarray*}}
\newcommand{\bdmn}{\begin{eqnarray}}
\newcommand{\edmn}{\end{eqnarray}}

\newcommand{\nn}{\nonumber}
\newcommand{\SSC}[1]{\section{#1}\setcounter{equation}{0}}

\newcounter{biblio}
\newenvironment{references}%
{\begin{list}{[\arabic{biblio}]}{\usecounter{biblio}%
\setlength{\leftmargin}{2.5em}\setlength{\rightmargin}{0pt}%
\setlength{\labelwidth}{2em}\setlength{\itemsep}{0pt}}}{\end{list}}
\newcommand{\References}%
{\vspace{2.8ex plus .3ex minus .3ex}%
\begin{center}{\bf References}\end{center}\begin{references}}

\usepackage{amssymb}
\usepackage{mathrsfs}
\newcommand{\bL}{{\mathbb{L}}}
\newcommand{\N}{{\mathbb{N}}}
\newcommand{\Z}{{\mathbb{Z}}}
\newcommand{\zd}{\Z^d}

\newcommand{\R}{{\mathbb{R}}}
\newcommand{\rd}{\R^d}

\newcommand{\ra }{\rightarrow }

\newcommand{\ov}{\overline}
\newcommand{\tl}{\widetilde}

\newcommand{\Llra}{\Longleftrightarrow }

\newcommand{\vvs}{\vspace{2ex}}
\newcommand{\vs}{\vspace{1ex}}

\newcommand{\lan}{\langle \:}
\newcommand{\ran}{\: \rangle}
\newcommand{\lef}{\left}
\newcommand{\rig}{\right}
\newcommand{\ri}{\right}
\newcommand{\st}{\stackrel}

\newcommand{\8}{\infty}

\newcommand{\dps}{\displaystyle}
\newcommand{\sub}{\subset}

\newcommand{\bsh}{\backslash}

\newcommand{\half}{\mbox{$\frac{1}{2}$}}

\newcommand{\inflim}{\mathop{\underline{\lim}}}

\renewcommand{\b}{\beta}

\newcommand{\del}{\delta}
\newcommand{\D}{\Delta}

\newcommand{\z}{\zeta}
\newcommand{\h}{\eta}
\newcommand{\tht}{\theta}

\newcommand{\lm}{\lambda}

\newcommand{\n}{\nu}

\newcommand{\rh}{\rho}

\newcommand{\s}{\sigma}

\renewcommand{\t}{\tau}

\newcommand{\W}{\Omega}

\newcommand{\cA }{{\cal A}}

\newcommand{\cF }{{\cal F}}

\newcommand{\cM }{{\cal M}}

\newcommand{\cP }{{\cal P}}

\newcommand{\cR }{{\cal R}}
\newcommand{\cS }{{\cal S}}

\makeatletter
\def\section{\@startsection{section}{1}{\z@}{-3.5ex plus -1ex minus 
 -.2ex}{2.3ex plus .2ex}{\bf}}
\def\subsection{\@startsection{subsection}{2}{\z@}{-3.25ex plus -1ex minus 
 -.2ex}{1.5ex plus .2ex}{\bf}}
\makeatother

\begin{document}
\bcenter

\large{\bf Localization for a Class of 
Linear Systems}\footnote{\today}

\vvs \normalsize

\vvs
Yukio Nagahata \\
Department of Mathematics, \\
Graduate School of Engineering Science \\
Osaka University,\\
Toyonaka 560-8531, Japan.\\
email: {\tt nagahata@sigmath.es.osaka-u.ac.jp}\\
URL: {\tt http://www.sigmath.osaka-u.ac.jp/}$\widetilde{}$ {\tt nagahata/}

\vvs
Nobuo Yoshida\footnote{
Supported in part by JSPS Grant-in-Aid for Scientific
Research, Kiban (C) 21540125} \\
Division of Mathematics \\
Graduate School of Science \\
Kyoto University,\\
Kyoto 606-8502, Japan.\\
email: {\tt nobuo@math.kyoto-u.ac.jp}\\
URL: {\tt http://www.math.kyoto-u.ac.jp/}$\widetilde{}$ {\tt nobuo/}

\ecenter
\begin{abstract}
We consider a class of continuous-time stochastic growth models on 
$d$-dimensional lattice with non-negative real numbers 
as possible values per site. 
The class contains examples such as binary contact path process 
and potlatch process.
We show the equivalence between the
slow population growth and localization property 
that the time integral of the 
replica overlap diverges. 
We also prove, under reasonable assumptions, a 
localization property in a stronger form that 
the spatial distribution of the population does not 
decay uniformly in space. 
\end{abstract}

\small
\noindent Abbreviated Title:  Localization for Linear Systems.
\tableofcontents

\normalsize
\SSC{Introduction}
We write $\N=\{0,1,2,...\}$, 
$\N^*=\{1,2,...\}$ and 
$\Z=\{ \pm x \; ; \; x \in \N \}$. For 
$x=(x_1,..,x_d) \in \rd$, $|x|$ stands for the $\ell^1$-norm: 
$|x|=\sum_{i=1}^d|x_i|$. For $\h=(\h_x)_{x \in \zd} \in \R^{\zd}$, 
$|\h |=\sum_{x \in \zd}|\h_x|$. 
Let 
$(\W, \cF, P)$ be 
a probability space. 
For events $A,B \sub \W$, $A \sub B$ a.s. means that 
$P(A \bsh B)=0$. Similarly, $A = B$ a.s. mean that 
$P(A \bsh B)=P(B \bsh A)=0$. By a {\it constant}, we always means a 
{\it non-random} constant. 

\vvs
We consider a class of continuous-time stochastic growth models on 
$d$-dimensional lattice $\zd$ with non-negative real numbers 
as possible values per site, so that the configuration at time $t$ 
can be written as $\h_t=(\h_{t,x})_{x \in \zd}$, $\h_{t,x} \ge 0$. 
We interpret the coordinate $\h_{t,x}$ as the ``population" at 
time-space $(t,x)$, though it need {\it not} be an integer. 
The class of growth models considered here is a reasonably 
ample subclass of the 
one considered in \cite[Chapter IX]{Lig85} as ``linear systems''. 
For example, it contains examples such as binary contact path process 
and potlatch process. The basic feature of the class is 
that the configurations are updated by applying 
the random linear transformation of the following form, when 
the Poisson clock rings at time-space $(t,z)$:
\bdnl{h_(x,t)intro}
\h_{t,x}=\lef\{ \barray{ll}
K_0\h_{t-,z} & \mbox{if $x=z$}, \\
\h_{t-,x}+K_{x-z}\h_{t-,z} & 
\mbox{if $x \neq z$},
\earray \ri.
\edn
where $K=(K_x)_{x \in \zd}$ is a random vector 
with non-negative entries, and independent copies of 
$K$ are used for each update (See section \ref{sec:model} for more 
detail). 
These models are known to exhibit, roughly speaking, 
the following phase transition \cite[Chapter IX, sections 3--5]{Lig85}:
\bds
\item[i)] 
If the dimension is high $d \ge 3$, {\it and} if the vector $K$ 
is not too random, then, with positive probability, 
the growth of the population is as fast 
as its expected value as time the $t$ tends to infinity, as such the 
{\it regular growth phase}.
\item[ii)] 
If the dimension is low $d=1,2$, {\it or} if the vector $K$ 
is random enough, then, almost surely, 
 the growth of the population strictly 
slower than its expected value as the time $t$ tends to 
infinity, as such the {\it slow growth phase}.
\eds
We denote the spatial distribution 
of the population by: 
\bdnl{rh}
\rh_{t,x}=\frac{\h_{t,x}}{|\h_t|}{\bf 1}_{\{ |\h_t|>0 \}}, \;
\; t >0, x \in \zd. 
\edn 
In our previous paper \cite{NY09}, we investigated the case (i) 
above and showed under some technical assumptions that the 
spatial distribution (\ref{rh}) obeys the central 
limit theorem. We also proved the delocalization property which 
says that the spatial distribution (\ref{rh}) decays uniformly 
in space like $t^{-d/2}$ as time $t$ tends to infinity.

In the present paper, 
we turn to the case (ii) above. We first prove the equivalence 
between the slow growth and a certain localization property in terms 
of the divergence of integrated replica overlap (\Thm{loc} below). 
We also show that, under reasonable assumptions, the localization 
occurs in stronger form that the spatial distribution (\ref{rh}) 
does not decay uniformly in space as time $t$ tends 
to infinity (\Thm{sloc} below). 
These, together with our previous work \cite{NY09}, verifies 
the delocalization/localization transition in correspondence with 
regular/slow growth transition for the class of model considered here.

It should be mentioned that the delocalization/localization transition 
in the same spirit has been discussed recently in various context, e.g.,  
\cite{CH02,CH06,CSY03,CY05,HY09,Sh09,Yo08a,Yo08b}. In particular, the 
last paper \cite{Yo08b} by the second author of the present article 
can be considered as the discrete-time counterpart 
of the present paper. Still, we believe it worth while verifying 
the delocalization/localization transition 
for the continuous-time growth models 
discussed here, in view of its classical importance of the model.
\subsection{The model} \label{sec:model}
We introduce a random vector 
$K=(K_x)_{x \in \zd }$ which is bounded and of 
finite range in the sense that
\bdnl{K_x}
0 \le  K_x \le b_K {\bf 1}_{\{ |x| \le r_K\}}\; \; 
\mbox{a.s. for some constants $b_K,r_K \in [0,\8)$.} 
\edn
Let $\t^{z,i}$, ($z \in \zd$, $i \in \N^*$) 
be i.i.d. mean-one exponential random variables 
and $T^{z,i}=\t^{z,1}+...+\t^{z,i}$. 
Let also $K^{z,i}=(K_x^{z,i})_{x \in \zd }$ 
($z \in \zd$, $i \in \N^*$) be i.i.d. 
random vectors with the same distributions as $K$, 
independent of $\{\t^{z,i} \}_{z \in \zd, i \in \N^*}$. 
Unless otherwise stated, we suppose for simplicity 
that the process $(\h_t)_{t \ge 0}$ 
starts from a single particle at the origin:
\bdnl{h_0}
\h_0=(\h_{0,x})_{x \in \zd},
\; \; \; 
\h_{0,x}=\lef\{ \barray{ll}
1 & \mbox{if $x=0$}, \\
0 & \mbox{if $x \neq 0$}.
\earray \rig.
\edn 
At time $t=T^{z,i}$, $\h_{t-}$ is replaced by $\h_t$, where 
\bdnl{h_(x,t)}
\h_{t,x}=\lef\{ \barray{ll}
K^{z,i}_0\h_{t-,z} & \mbox{if $x=z$}, \\
\h_{t-,x}+K^{z,i}_{x-z}\h_{t-,z} & 
\mbox{if $x \neq z$}.
\earray \ri.
\edn
A formal construction of the process $(\h_t)_{t \ge 0}$ can be 
given as a special case of  \cite[page 427, Theorem 1.14]{Lig85} 
via Hille-Yosida theory. In section \ref{p|ovh|=M}, 
we will also give an alternative construction of the process 
in terms of a stochastic differential equation.

To exclude uninteresting cases from the 
viewpoint of this article, we also assume that
\bdmn 
& & \mbox{the set $\{x \in \zd\; ;\; E[K_x]\neq 0\}$ contains a 
linear basis of $\rd$,} \label{true_d} \\
& & P(|K|= 1)<1. \label{Kneq1}
\edmn
The first assumption (\ref{true_d}) makes the model 
``truly $d$-dimensional". The reason for the second assumption 
(\ref{Kneq1}) is to exclude the case $|\h_t| \equiv 1$ a.s.

Here are some typical examples which fall into the above set-up:

\vs
\noindent $\bullet$
{\bf The binary contact path process (BCPP):}
The binary contact path process (BCPP), originally introduced by 
D. Griffeath \cite{Gri83} is a special case the model, where 
\bdnl{binK} 
K= \lef\{ \barray{ll}
\lef(\del_{x,0}+\del_{x,e} \rig)_{x \in \zd} & 
\mbox{with probability ${\lm \over 2d\lm +1}$, for each $2d$ neighbor 
$e$ of 0} \\
0 & \mbox{with probability ${1 \over 2d\lm +1}$}. 
\earray  \rig. \edn
The process is interpreted as the spread of an infection, 
with $\h_{t,x}$ infected individuals at time $t$ at the site $x$. 
The first line of (\ref{binK}) says that, 
with probability ${\lm \over 2d\lm +1}$ for each $|e|=1$, 
all the infected individuals at site $x-e$ are 
duplicated and added to those on the site $x$.
On the other hand, the second line of (\ref{binK}) says that, 
all the infected individuals at a site become healthy with 
probability ${1 \over 2d\lm +1}$. 
A motivation to study the BCPP comes from the 
fact that the projected process
$$
\lef( \h_{t,x} \wedge 1\ri)_{x \in \zd},\; \; \; t \ge 0
$$
is the basic contact process \cite{Gri83}.

\vs
\noindent $\bullet$
{\bf The potlatch process:}
The potlatch process discussed in e.g. 
\cite{HL81} and \cite[Chapter IX]{Lig85} 
is also a special case of the above set-up, 
in which
\bdnl{binP/S} 
K_x=W k_x,\; \; x \in \zd.
\edn
Here, $k=(k_x)_{x \in \zd} \in [0,\8)^{\zd}$ is a 
non-random vector and $W$ is a non-negative, bounded, mean-one 
random variable such that $P(W = 1)<1$ 
(so that the notation $k$ here is consistent with the definition 
(\ref{kp12}) below). The potlatch process 
was first introduced in 
\cite{Spi81} for the case $W \equiv 1$ and discussed further in 
\cite{LS81}. It was in \cite{HL81} where case with $W \not\equiv 1$ was 
introduced and discussed. 
Note that we do not restrict 
ourselves to the case $|k|=1$ unlike in 
\cite{HL81} and \cite[Chapter IX]{Lig85}. 
\subsection{The regular and slow growth phases}

We now recall the following facts and notion from 
\cite[page 433, Theorems 2.2 and 2.3]{Lig85}, 
although our terminologies are somewhat different from 
the ones in \cite{Lig85}.
  Let $\cF_t$ be the 
$\s$-field generated by $\h_s$, $s \le t$.
\Lemma{0,1}
We set
\bdmn
k&=& (k_x)_{x \in \zd}=(E[K_x])_{x \in \zd} \; \; \; \label{kp12} \\
\ov{\h}_t &=&(e^{-(|k|-1)t }\h_{t,x})_{x \in \zd}. \label{ovh_t}
\edmn
Then,
\bds
\item[a)]
$(|\ov{\h}_t|, \cF_t)_{t \ge 0}$ is a martingale, 
and therefore, the following limit exists a.s. 
\bdnl{ovn_8}
|\ov{\h}_\8| =\lim_{t \ra \8}|\ov{\h}_t|.
\edn
\item[b)] Either 
\bdnl{0,1}
E[|\ov{\h}_\8|]=1\; \; \mbox{or}\; \; 0.
\edn
Moreover, $E[|\ov{\h}_\8|]=1$ if and only if the limit (\ref{ovn_8}) 
is convergent in $\bL^1 (P)$. 
\eds
\end{lemma}
We will refer to the former case of (\ref{0,1})
as {\it regular growth phase} and the latter as 
{\it slow growth phase}. 

The regular growth 
means that, at least with positive probability, 
the growth of the ``total number" $|\h_t|$ 
of the particles is of the same order as its expectation 
$e^{(|k|-1)t}|\h_0|$.
On the other hand, the slow growth means that, almost surely, 
 the growth of $|\h_t|$ is slower than its expectation. 

Since we are mainly interested in the slow growth phase in this paper, 
we  now present sufficient conditions for the slow growth.
\Proposition{SG}
\bds 
\item[a)]
For $d=1,2$, ${\dps |\ov{\h}_\8|=0}$ a.s. 
In particular for $d=1$, there exists a constant $c>0$ such that 
\bdnl{SG1}
|\ov{\h}_t| =O (e^{-ct}),\; \; \; \mbox{as $t \ra \8$, a.s.}
\edn 
\item[b)] 
For any $d \ge 1$, suppose that 
\bdnl{SGlog}
\sum_{x \in \zd}E\lef[ K_{x}\ln K_{x}\ri] >|k|-1
\edn
Then, again, 
there exists a constant $c>0$ such that (\ref{SG1}) holds.
\eds
\end{proposition}
Proof: Except for (\ref{SG1}), these sufficient conditions are 
presented in 
\cite[Chapter IX, sections 4--5]{Lig85}. The exponential decay 
(\ref{SG1}) follows from similar arguments as in discrete-time 
models discussed in \cite[Theorems 3.1.1 and 3.2.1]{Yo08a}. 
\hfill $\Box$

\vvs
\noindent {\bf Remarks:}
{\bf 1)}
For BCPP, (\ref{SGlog}) is equivalent to $\lm <(2d)^{-1}$, 
in which case it is known that $|\h_t| \equiv 0$ 
for large enough $t$'s a.s. 
\cite[Example 4.3.(c) on page 33, together with 
Theorem 1.10 (a) on page 267]{Lig85}.
Thus, \Prop{SG}(b) applies only in a trivial manner for BCPP. 
In fact, we do not know if there is a value $\lm$ for 
which BCPP with $d \ge 3$ is in slow growth phase, 
without getting extinct a.s. 
For potlatch  process, 
$$
(\ref{SGlog})\; \; \Llra \; \; 
E[W \ln W] >{|k|-1-\sum_x k_x \ln k_x \over |k|}.
$$
Thus, (\ref{SGlog}) and hence  (\ref{SG1}) is true if 
$W$ is ``random enough". \\
\noindent {\bf 2)} 
A sufficient condition for the regular growth phase 
will be given by (\ref{sloc_sharp}) below.
\subsection{Results}
Recall that we have defined the spatial distribution 
of the population by (\ref{rh}).
Interesting objects related to the
density would be 
\bdnl{rh^*} \rh^*_t=\max_{x \in \zd}\rh_{t,x}, \;
\; \mbox{and}\; \; \cR_t=\sum_{x \in \zd}\rh_{t,x}^{2}. 
\edn
$\rh^*_t$ is the density at the most populated site, while $\cR_t$
is the probability that a given pair of particles at time $t$ are
at the same site. We call $\cR_t$ the {\it replica overlap}, in
analogy with the spin glass theory. Clearly, $(\rh^*_t)^{2} \le \cR_t
\le \rh^*_t$. These quantities convey information on
localization/delocalization of the particles. Roughly speaking,
large values of $\rh^*_t$ or $\cR_t$ indicate that the most of
the particles are concentrated on small number of ``favorite
sites" ({\it localization}), whereas small values of them imply
that the particles are spread out over a large number of sites ({\it
delocalization}).

We first show that
the regular and slow growth are characterized, respectively  by 
convergence (delocalization) and divergence (localization) 
of the integrated replica overlap: $\int^\8_0\cR_s ds$. 
\Theorem{loc}
\bds
\item[a)]
Suppose that $P(|\ov{\h}_\8|>0)>0$. Then, 
$$
\int^\8_0\cR_s ds <\8\; \; \mbox{a.s.}
$$
\item[b)] 
Suppose on the contrary that $P(|\ov{\h}_\8|=0)=1$. Then, 
\bdnl{S=L}
\{ \; \mbox{survival}\; \} 
= 
\lef\{ \; \int^\8_0\cR_s ds=\8 \; \ri\},
\; \; \; \mbox{a.s.}
\edn
where $\{ \mbox{survival}\}=\{ |\h_t| \neq 0 \; \mbox{for all $t \ge 0$} \}$.
Moreover, there exists a constant $c>0$ such that 
\bdnl{M<}
|\ov{\h}_t| \le  \exp \lef( -c
\int^t_0\cR_s ds\ri) 
\; \; \mbox{for all large enough $t$'s, a.s.}
\edn
\eds
\end{theorem}
Results of this type are fundamental in analyzing a certain 
class of spatial random growth models, 
such as directed polymers in random environment 
\cite{CH02,CH06,CSY03,CY05}, 
linear stochastic evolutions \cite{Yo08b}, branching 
random walks and Brownian motions in random environment \cite{HY09,Sh09}.
Until quite recently, however, this type of results were available only 
when no extinction at finite time is allowed, i.e.,
  $|\h_t|>0$ for all $t \ge 0$, e.g., 
\cite{CH02,CH06,CSY03,CY05,HY09,Sh09}.  
In fact, the proof there relies on the analysis of the 
supermartingale $\ln |\bar{\h}_t|$, which is not even 
defined if extinction at finite time is possible. 
To overcome this problem, we will adapt 
a more general approach introduced in \cite{Yo08b}. 

\vvs
Next, we present a result (\Thm{sloc} below) 
which says that, under reasonable assumptions, 
we can strengthen the localization property
$$
\int_0^\8 \cR_s ds=\8 
$$
 in (\ref{S=L}) to:
$$
\int_0^\8 {\bf 1}\{\cR_s \ge c \} ds=\8,
$$
where $c>0$ is a constant. 
To state the theorem, we define
\bdnl{beta_x,y}
\b_{x,y}=E[(K-\del_0)_x(K-\del_0)_y], \; \; \; 
x,y \in \zd.
\edn
We also introduce:
\bdnl{G(x)}
G(x)=\int^\8_0P_S^0(S_t=x)dt,
\edn
where $((S_t)_{t \ge 0}, P_S^x)$ is the continuous-time random walk on 
$\zd$ starting from $x \in \zd$, with the generator 
\bdnl{L_S}
L_Sf (x)=\half \sum_{y \in \zd}\lef( k_{x-y}+k_{y-x} \ri)
\lef( f(y)-f(x)\rig), \; \; \; \mbox{cf. (\ref{kp12}).}
\edn
\Theorem{sloc}
Referring to (\ref{beta_x,y})--(\ref{G(x)}), suppose either of 
\bds
\item[a)]
$d=1,2$.
\item[b)] $d \ge 3$, $P (|\ov{\h}_\8| =0)=1$ and 
\bdnl{cond:sloc}
\sum_{x,y \in \zd}G(x-y)\b_{x,y}>2.
\edn
\eds
Then there exist a constant $c \in (0,1]$ 
such that
\bdnl{sloc}
\{ \;\mbox{survival}\; \} 
= \lef\{
\int_0^\8 {\bf 1}\{\cR_s \ge c \} ds=\8 \; \ri\}
 \; \; a.s. 
\edn
\end{theorem}
Our proof of \Thm{sloc} is based on the idea of P. Carmona and 
Y. Hu in \cite{CH02,CH06}, where they prove similar results 
for directed polymers in random environment. Although the arguments in 
\cite{CH02,CH06} are rather complicated and uses special 
structure of the model, it was possible to  
extract the main idea from \cite{CH02,CH06} in a way applicable 
to our setting. Also, we could considerably reduce the technical 
complexity in the argument as compared with \cite{CH02,CH06}.

\vvs
\noindent {\bf Remarks}:
\noindent {\bf 1)}
We prove (\ref{sloc}) by way of the following stronger estimate:
\bdnl{sloc2}
\{ \;\mbox{survival}\; \} 
\sub \lef\{ \; 
\inflim_{t \nearrow \8} {\int_0^t\cR_s^{3/2}ds 
\over \int_0^t\cR_sds} \ge c_1\; \ri\} \; \; \mbox{a.s}. 
\edn
for some constant  $c_1>0$. The inequality $r^{3/2} 
\le {\bf 1}\{ r \ge c\}+\sqrt{c}r$ for $r,c \in [0,1]$ can be used 
to conclude (\ref{sloc}) from (\ref{sloc2}). \\
\noindent {\bf 2)}
We note that $P (|\ov{\h}_\8|>0)>0$ if 
\bdnl{sloc_sharp}
d \ge 3 \; \; \mbox{and}\; \; \sum_{x,y \in \zd}G(x-y)\b_{x,y}<2.
\edn
This, together with \Thm{loc}(a), shows that 
the condition (\ref{cond:sloc}) is necessary, up to the equality, 
for (\ref{sloc}) to be true whenever 
survival occurs with positive probability. 
We see that (\ref{sloc_sharp}) implies $P (|\ov{\h}_\8|>0)>0$ 
via the same line of argument as in 
\cite[page 464, Theorem 6.16]{Lig85}, where the special case of 
the potlatch process is discussed.
 We consider the {\it dual process}
$\z_t \in [0,\8)^{\zd}, \; \; t \ge 0$ which evolves in the 
same way as $(\h_t)_{t \ge 0}$ except that 
(\ref{h_(x,t)intro}) is replaced by its transpose:
\bdnl{h_(x,t)*}
\z_{t,x}=\lef\{ \barray{ll}
\sum_{y \in \zd}K_{y-x}\z_{t-,y} & \mbox{if $x=z$}, \\
\z_{t-,x} & \mbox{if $x \neq z$}.
\earray \ri.
\edn 
By \cite[page 445, Theorem 3.12]{Lig85}, a sufficient condition 
for $P (|\ov{\h}_\8|>0)>0$ is that there exists a 
function $h:\zd \ra (0,\8)$ 
such that $\lim_{|x| \ra \8}h(x)= 1$ and that 
\bdnl{sum_yq(x,y)h(y)=0}
\sum_yq(x,y)h(y)=0,\; \; \; x \in \zd.
\edn
Here, $q(x,y)$ is the matrix given by 
\cite[page 445, (3.8)--(3.9)]{Lig85} for the dual process. In our setting, 
it is computed as:
$$
q(x,y)=k_{x-y}+k_{y-x}-2|k|\del_{x,y}+\del_{0,x}\sum_z\b_{z,z+y},
$$
so that (\ref{sum_yq(x,y)h(y)=0}) becomes:
$$
(L_S h)(x)+\half \del_{0,x}\sum_{y,z}h(y-z)\b_{y,z}=0, \; \; \; x\in \zd, 
\; \; \; \mbox{cf. (\ref{L_S}).}
$$
Under the assumption (\ref{sloc_sharp}), 
a choice of such function $h$ is given by 
$h=1+cG$, where 
$$
c={E[(|K|-1)^2] \over 1-\half \sum_{x,y \in \zd}G(x-y)\b_{x,y}}.
$$
\noindent {\bf 3)}
Let $\pi_d$ be the return probability for the simple random 
walk on $\zd$. Also, let $\lan \cdot, \; \cdot \ran$ and $*$ be the 
inner product of $\ell^2 (\zd)$ and the discrete convolution respectively.
We then have that
\bdnl{cond:sloc_exam}
\mbox{(\ref{cond:sloc}) }\; \; \Llra \; \; 
\lef\{ \barray{ll}
\lm <{1 \over 2d (1-2\pi_d)} & \mbox{for BCPP},\\
E[W^2] >{(2|k|-1)G(0)\over \lan G*k, k \ran}
& \mbox{for the potlatch process.} 
\earray \rig.
\edn
For BCPP, (\ref{cond:sloc_exam}) can be seen from that 
(cf. \cite[page 965]{NY09})
$$
\b_{x,y}={{\bf 1}\{ x=0\}+\lm {\bf 1}\{ |x|=1\} \over 2d\lm +1}\del_{x,y}, 
\; \; \mbox{and}\; \; 
G(0)={2d\lm +1 \over 2d\lm}{1 \over 1-\pi_d}.
$$
To see (\ref{cond:sloc_exam}) for the potlatch process, 
we note that $\half (k+ \check{k})*G =|k|G-\del_0$, 
with $\check{k}_x=k_{-x}$ and that
$$
\b_{x,y}=E[W^2]k_xk_y-k_x \del_{y,0}-k_y \del_{x,0}+\del_{x,0}\del_{y,0}.
$$ 
Thus, 
\bdnn
\sum_{x,y \in \zd}G(x-y)\b_{x,y}
& = & E[W^2]\lan G*k, k \ran -\lan G, k+ \check{k} \ran +G(0) \\
& = & E[W^2]\lan G*k, k \ran +2-(2|k|-1)G(0),
\ednn
from which (\ref{cond:sloc_exam}) for the potlatch process follows. 
\subsection{SDE description of the process} \label{p|ovh|=M}
We now give an alternative description of the process in terms of a 
stochastic differential equation (SDE).
We introduce random measures on $[0,\8) \times [0,\8)^{\zd}$ by 
\bdnl{N^z}
N^z( dsd\xi)
=\sum_{i \ge 1}{\bf 1}\{ (T^{z,i}, K^{z,i}) \in dsd\xi \}, 
\; \; \; N^z_t(dsd\xi)={\bf 1}_{\{s \le t\}}N^z(dsd\xi).
\edn
Then, $N^z$, $z \in \zd$ are independent Poisson random measures on 
$[0,\8) \times [0,\8)^{\zd}$ with the intensity 
$$
ds \times P (K \in \cdot).
$$ 
The precise definition of the process $(\h_t)_{t \ge 0}$ is then  given by 
the following stochastic differential equation:
\bdnl{sde}
\h_{t,x}=\h_{0,x}
+\sum_{z \in \zd}\int N^z_t(dsd\xi)
\lef( \xi_{x-z}-\del_{x,z}\ri)\h_{s-,z}.
\edn
By (\ref{K_x}), it is standard to see that (\ref{sde}) defines  
a unique process $\h_t=(\h_{t,x})$, ($t \ge 0$) 
and that $(\h_t)$ is Markovian. 
\SSC{Proofs}
It is convenient to introduce the following notation:
\bdmn
& & \n =P(K \in \cdot) \in \cP ([0,\8)^{\zd}), \; \; \mbox{the law of $K$}.
\label{nu} \\
& & \tl{N}^z(dsd\xi)=N^z(dsd\xi)-ds\n (d\xi), \; \; \; 
\tl{N}^z_t(dsd\xi)={\bf 1}_{s \le t}\tl{N}^z(dsd\xi).
\label{tlN}
\edmn
\subsection{Proof of \Thm{loc}}
The proof of \Thm{loc} is based on the following
\Lemma{SS=L}
\bdnl{loc}
\{ \;|\ov{\h}_\8|=0, \; \; \mbox{survival}\; \} 
= 
\lef\{ \; \int^\8_0\cR_s ds=\8 \; \ri\},
\; \; \; \mbox{a.s.}
\edn
Moreover, there exists a constant $c>0$ such that 
(\ref{M<}) holds 
a.s. on the event $\lef\{ \; \int^\8_0\cR_s ds=\8 \; \ri\}$.
\end{lemma}
Proof: 
We see from (\ref{sde}) that
\bdnn
|\ov{\h}_t|
&=& |\ov{\h}_0|
+\sum_z\int \tl{N}^z_t (ds d\xi)|\ov{\h}_{s-}|(|\xi|-1)\rh_{s-,z}
\; \; \;  \mbox{(cf. (\ref{tlN}))}\\
&=& |\ov{\h}_0|+\int^t_0|\ov{\h}_{s-}| dM_{s}
\ednn
where
$$
M_{t} = \sum_z\int \tl{N}^z_t (ds d\xi)(|\xi|-1)\rh_{s-,z}.
$$
Then, by the Dol\'eans-Dale exponential formula 
(e.g., \cite[page 248, 9.39]{HWY92}), 
$$
|\ov{\h}_t| =
\exp \lef( M_t \ri)D_t,
$$
where
$$
D_t=\prod_{s \le t} \lef( 1+\D M_{s} \ri)\exp \lef( -\D M_{s}\ri), 
\; \; \mbox{with}\; \; \D M_t=M_t-M_{t-}.
$$
Note also the predictable quadratic variation of $M_\cdot$ is given 
by 
\bds \item[1)] \hspace{1cm}${\dps 
\lan M \ran_t =E[(|K|-1)^2]\int^t_0\cR_s ds.
}$ \eds
Since $-1 \le \D M_t \le b_K-1 <\8$, we have that 
(See e.g.\cite[page 222, 8.32]{HWY92})
\bds \item[2)] \hspace{1cm}${\dps 
\{ \; \lan M \ran_\8<\8 \; \}
\sub \{ [M]_\8 <\8, \; \mbox{$M_{t}$ converges as $t \nearrow \8$}\}
\; \; }$ a.s.
\item[3)] \hspace{1cm}${\dps 
\{ \; \lan M \ran_\8=\8 \; \}  
\sub
\lef\{ \lim_{t \ra \8} {\lan M \ran_t \over [M]_t}=1, 
\; \; \lim_{t \ra \8} {M_t \over \lan M \ran_t}=0\; \ri\} \; \; }$ a.s.
\eds
where 
$$
[M]_t=\sum_{s  \le t}(\D M_{s})^2
$$
We start with the ``$\supset$'' part of (\ref{loc}): 
Note that 
$(1+u)e^{-u} \le e^{-c_1u^2}$ for $-1 \le u \le b_K-1$, where  
$c_1>0$ is a constant. 
We suppose that $\int^\8_0\cR_s ds=\8$, or equivalently that, 
$\lan M \ran_\8 =\8 $.
Then, for large $t$, 
$$
\exp \lef( M_t \ri)D_t \le \exp \lef( M_t -c_1[M]_t \ri)
\st{\mbox{\scriptsize (3)}}{\le}
\exp \lef( -{c_1 \over 2}\lan M \ran_t \ri)
\st{\mbox{\scriptsize (1)}}{\le}
 \exp \lef( -c_2 \int^t_0\cR_s ds \ri)
$$
This shows that $\int^\8_0\cR_s ds=\8$ implies $|\h_\8|=0$, 
together with the bound (\ref{M<}). \\
We now turn to the ``$\sub$'' part of (\ref{loc}):  
We need to prove that
\bds
\item[4)]
$\{ \; \int^\8_0\cR_s ds<\8 \; \; \; \mbox{survival}\} 
\st{\rm a.s.}{\sub} \{ \;|\ov{\h}_\8|>0\} $.
\eds 
We have 
\bds
\item[5)]
$\{ \; \int^\8_0\cR_s ds<\8 \}  
\st{\mbox{\scriptsize (1)--(2)}}{\sub} 
\{ \mbox{$M_{t}$ converges as $t \nearrow \8$}\}\; \; $ a.s.
\eds
On the other hand, 
$$
\sum_{s \le t}
\lef| \lef( 1+\D M_{s} \ri)\exp \lef( -\D M_{s} \ri)-1 \ri|
\le {e \over 2} [M ]_t,
$$
since $|(1+u)e^{-u}-1| \le eu^2/2$ for $u \ge -1$. Thus, 
\bds
\item[6)]
$\{ \; \int^\8_0\cR_s ds<\8, \; \; \; \mbox{survival}\} \sub
\lef\{ \mbox{$D_t$ converges to a positive limit 
as $t \nearrow \8$}\ri\}\; \; $ a.s.
\eds 
We now obtain (4) by (5)--(6).
\hfill $\Box$

\vvs
\noindent {\it Proof of \Thm{loc}}:
a):
If $P(|\ov{\h}_\8|>0)>0$, then, 
$$
\lef\{ \mbox{survival} \rig\}=
\{ |\ov{\h}_\8|>0\}\; \; \;\mbox{a.s.}
$$
This can be seen easily by the argument in 
\cite[page 701, proof of Proposition]{Gri83}. We see from 
this and (\ref{loc}) that $\int^\8_0\cR_s ds <\8$ a.s. 
on the event of survival, 
while $\int^\8_0\cR_s ds <\8$ is obvious outside the event of survival. \\
b): This follows from \Lem{SS=L}
\hfill $\Box$

\subsection{Proof of \Thm{sloc}} \label{p:sloc1}
Let $p$ be a transition function of a symmetric discrete-time 
random walk 
defined by 
\[
p(x) = \lef\{ \barray{ll} 
\frac{\dps k_x + k_{-x}}{\dps 2(|k|-k_0)} & \mbox{if $x \ne 0$,} \\
0 & \mbox{if $x =0$.}
\earray \rig.
\]
and $p_n$ be the $n$-step transition function. 
We set 
\[
g_n(x) = \del_{x,0}+\sum_{k=1}^n p_k(x).
\]
\Lemma{g_n}
Under the assumptions of \Thm{sloc}, there exists $n$ such that
\bdnl{g_n}
\sum_{x,y}g_n(x-y)\b_{x,y}>2(|k|-k_0).
\edn
\end{lemma}
Proof: 
Since the discrete-time 
random walk with the transition probability $p$ is 
the jump chain of the continuous-time random walk 
$((S_t)_{t \ge 0}, P_S^x)$ with the generator (\ref{L_S}), 
we have that
\bds \item[1)] \hspace{1cm} ${\dps 
\lim_{n \ra \8}g_n(x)=(|k|-k_0)G(x)\; \; \;}$
for all $x \in \zd$.
\eds 
For $d \ge 3$, $G(x)<\8$ for any $x \in \zd$ and 
$\b_{x,y}\neq 0$ only when $|x|, |y| \le r_K$, we see from (1) that
$$
\lim_{n \ra \8}\sum_{x,y}g_n(x-y)\b_{x,y}
=(|k|-k_0)\sum_{x,y}G(x-y)\b_{x,y}.
$$
Thus, (\ref{g_n}) holds for all large enough $n$'s. \\
To show (\ref{g_n}) for $d=1,2$, we will prove that
$$
\lim_{n \ra \8}\sum_{x,y}g_{2n-1}(x-y)\b_{x,y}=\8.
$$
For $f \in \ell^1 (\zd)$, we denote 
its Fourier transform by 
$$
\widehat{f}(\tht)=\sum_{x \in \zd}f(x)\exp ({\bf i} x \cdot \tht),
 \; \; \; \tht \in I \st{\rm def}{=}[-\pi,\pi]^d.
$$
We then have that
$$
g_{2n-1}(x)={1 \over (2\pi)^d}\int_I
{1 -\widehat{p}(\tht)^{2n} \over 1 -\widehat{p}(\tht)}
\exp ({\bf i} x \cdot \tht)d\tht
$$
and hence that
\bdnn
\sum_{x,y}g_{2n-1}(x-y)\b_{x,y} 
&=&{1 \over (2\pi)^d}\int_I
{1 -\widehat{p}(\tht)^{2n} \over 1 -\widehat{p}(\tht)}
\sum_{x,y}\exp ({\bf i} (x-y) \cdot \tht)
E[(K-\del_0)_x (K-\del_0)_y]d\tht \\
&=&{1 \over (2\pi)^d}\int_I
{1 -\widehat{p}(\tht)^{2n} \over 1 -\widehat{p}(\tht)}
E[|\widehat{K}(\tht)-1|^2]d\tht .
\ednn
Since $p(\cdot)$ is even, 
we see that $ \widehat{p}(\tht) \in [-1,1]$ for all $\tht \in I$.
Also, by (\ref{true_d}), 
there exist constants $c_i>0$ ($i=1,2,3$) such that 
$$
0 \le 1-c_1|\tht|^2 \le \widehat{p}(\tht) \le 1-c_2|\tht|^2
\; \; \; \mbox{for $|\tht| \le c_3$.}
$$ 
These imply that
$$
\inflim_{n \ra \8}
\sum_{x,y}g_{2n-1}(x-y)\b_{x,y}
\ge {1 \over (2\pi)^dc_1}\int_{|\tht| \le c_2}
{ E[|\widehat{K}(\tht)-1|^2]\over |\tht|^2}d\tht.
$$
The integral on the right-hand-side diverges if $d \ge 2$, since  
$$
E[|\widehat{K}(0)-1|^2]=E[(|K|-1)^2] \neq 0.
$$
\hfill $\Box$

\vvs
We take an $n$ in \Lem{g_n} and fix it.  
We then set 
\bdnl{g&S}
g=g_n\; \; \mbox{and}\; \; \cS_t =
\lan g*\rh_t, \rh_t \ran,
\edn
where the bracket $\lan \cdot, \; \cdot \ran$ and $*$ stand for the 
inner product of $\ell^2 (\zd)$ and the discrete convolution respectively. 
In what follows, we will often use the Hausdorff-Young inequality:
\bdnl{HYineq}
|(f*h)^2|^{1/2}\le |f||h^2|^{1/2}\; \; \; 
f \in \ell^1(\zd), \; h \in \ell^2(\zd).
\edn
For example, we have that
\bdnl{|X_t|}
0 \le \cS_t \st{\mbox{\scriptsize Schwarz}}{\le}
  |(g*\rh_t)^2|^{1/2}|(\rh_t)^2|^{1/2}
\st{\mbox{\scriptsize (\ref{HYineq})}}{\le}
|g||(\rh_t)^2| = |g| \cR_t< \infty.
\edn
The proof of \Thm{sloc} is based on the following
\Lemma{cS=cM+cA}
Let 
$$
\cS_t = \cS_0+\cM_t + \cA_t
$$
be the Doob decomposition, where $\cM_\cdot$ and $\cA_\cdot$ 
are a martingale and a predictable process, respectively. Then,
\bds
\item[a)] 
There is constants $c_1,c_2 \in (0,\8)$ such that
\bdnl{cA>}
\cA_t \ge \int_0^t \lef( c_1 \cR_s -c_2 \cR_s^{3/2} \ri) ds
\edn
\item[b)]
\bdnl{M=o(intR)}
\lef\{ \int^\8_0 \cR_s ds =\8 \ri\} \sub
\lef\{ \lim_{t \ra \8}{\cM_t \over \int^t_0 \cR_s ds}=0 \ri\}
\; \; \; \mbox{a.s.}
\edn
\eds
\end{lemma}
{\it Proof of \Thm{sloc}:}
By \Thm{loc} and the remark after \Thm{sloc}, it is enough to 
prove that
\bds \item[1)] \hspace{1cm}
${\dps   
\inflim_{t \nearrow \8} 
{\int_0^t\cR_s^{3/2}ds \over \int_0^t\cR_sds} \ge c \; \; }$ a.s. 
on ${\dps D\st{\rm def}{=}
\lef\{\int_0^\infty \cR_t dt = \infty \ri\}}$
\eds
for a positive constant $c$.
It follows from  (\ref{|X_t|}) and (\ref{M=o(intR)}) that 
$$
\lim_{t \ra \8}{\cA_t \over \int^t_0 \cR_s ds}=0
\; \; \; \mbox{a.s. on $D$}
$$
and hence from (\ref{cA>}) that
$$
\inflim_{t \ra \8}{\int^t_0 \cR_s^{3/2} ds 
\over \int^t_0 \cR_s ds} \ge {c_1 \over c_2}
\; \; \; \mbox{a.s. on $D$}.
$$
This proves (1) and hence \Thm{sloc}.
\hfill $\Box$
\subsection{Proof of \Lem{cS=cM+cA}}
\noindent {\it Proof of part (a)}:
To make the expressions below easier to read, we introduce the 
following shorthand notation:
\bdnn
J_{t,x,z}(\xi)
&=&\rh_{t,x}+(\xi - \del_{0})_{x-z} \rh_{t,z}, \\
\ov{J}_{t,x,z}(\xi)
&=&
{\h_{t,x}+(\xi - \del_{0})_{x-z} \h_{t,z} \over 
|\h_t|+(|\xi|-1)\h_{t,z} }=
{J_{t,x,z}(\xi)\over 1+(|\xi|-1)\rh_{t,z} }.
\ednn
We then rewrite $\cS_t$ as:
\begin{eqnarray*}
\cS_t &=& \cS_0 + \sum_z \int N_t^z(du d\xi)
\sum_{x,y}g(x-y)
\lef(
  \ov{J}_{u-,x,z}(\xi) \ov{J}_{u-,y,z}(\xi)-\rho_{u-,x} \rho_{u-,y}
\ri) \\
&=& \cS_0 + \cM_t + \cA_t
\end{eqnarray*}
where $\cA_t=\int^t_0A_sds$ has been defined by 
$$
A_s = \sum_{x,y,z}g(x-y)\int \nu(d\xi)
\lef(
  \ov{J}_{s,x,z}(\xi) \ov{J}_{s,y,z}(\xi)-\rho_{s,x} \rho_{s,y}
\ri)
$$
To bound $A_s$ from below, we note that 
$(1+x)^{-2} \ge 1 - 2x$ for $x \ge -1$. Then,
\bdmn \lefteqn{
 \ov{J}_{s,x,z}(\xi) \ov{J}_{s,y,z}(\xi)-\rho_{s,x} \rho_{s,y} }\nn \\
& \ge & J_{t,x,z}(\xi)J_{t,y,z}(\xi)
-2(|\xi|-1)\rh_{s,z}J_{t,x,z}(\xi)J_{t,y,z}(\xi)-\rho_{s,x} \rho_{s,y} \nn \\
& = & U_{s,x,y,z}(\xi)-2V_{s,x,y,z}(\xi)-2W_{s,x,y,z}(\xi),
\label{JJ-rr}
\edmn
where 
\bdmn
U_{s,x,y,z}(\xi)
&=& J_{s,x,z}(\xi)J_{s,y,z}(\xi)-\rho_{s,x} \rho_{s,y} \label{Utxyz}  \\
V_{s,x,y,z}(\xi) &=& (|\xi|-1) U_{s,x,y,z}(\xi)\rh_{s,z}
\label{Vtxyz} \\
W_{s,x,y,z}(\xi)
&= & (|\xi|-1)\rh_{s,x} \rh_{s,y} \rh_{s,z}.
\label{Wtxyz}
\edmn
We will see that
\bdnl{intV}
\sum_{x,y,z}g(x-y)\int V_{s,x,y,z}(\xi)\nu(d\xi)
\le c\cR_s^{3/2}.
\edn 
Here and in what follows, $c$ denotes a multiplicative constant, 
which does not depends on time variable $s$ and space variables 
$x,y,...$.
To prove (\ref{intV}), we can bound the factor $|\xi|-1$ by a 
constant. We write 
\bdnl{Uexp}
U_{s,x,y,z}(\xi)=
(\xi - \del_{0})_{y-z}\rho_{s,x} \rho_{s,z}
+(\xi - \del_{0})_{x-z}\rho_{s,y} \rho_{s,z} 
 +(\xi - \del_{0})_{x-z}(\xi - \del_{0})_{y-z}\rho_{s,z}^2 
\edn
We look at the contribution from the second term 
on the right-hand-side of (\ref{Uexp}) to the 
left-hand-side of (\ref{intV}). 
\bdnn
\sum_{x,y,z}g(x-y)(\xi - \delta_0)_{x-z} \rho_{s,z}^2 \rho_{s,y}
&=&\lan g * \rh_s, (\xi - \delta_0) * \rh_s^2 \ran \\
&\le &|(g * \rh_s)^2|^{1/2}
|((\xi - \delta_0) * \rh_s^2)^2|^{1/2} \\
& \le & |g|\cR_s^{1/2} |(\xi - \delta_0)^2|^{1/2}
|\rh_s^2| \le c\cR_s^{3/2}
\ednn
Contributions 
from the other two terms on the right-hand-side of (\ref{Uexp}) 
can be bounded similarly. Hence we get (\ref{intV}). 

On the other hand, 
\bdmn
\lefteqn{\sum_{x,y,z}g(x-y)\int U_{s,x,y,z}d\nu} \nn \\
&=& 
\sum_{x,y,z}g(x-y)\Big(
(k - \del_{0})_{y-z}\rho_{s,x} \rho_{s,z} 
 +(k - \del_{0})_{x-z}\rho_{s,y} \rho_{s,z} 
+\b_{x-z,y-z}\rho_{s,z}^2 \Big) \nn \\
&=& 
\lan g*(k - \del_{0})*\rh_s, \rh_s \ran 
+\lan g*(\check{k} - \del_{0})*\rh_s, \rh_s \ran 
+\sum_{x,y}g(x-y)\b_{x,y}\cR_s,  \label{intU}
\edmn
where 
$\check{k}_x=k_{-x}$. 
Also, 
\bdnl{intW}
\sum_{x,y,z}g(x-y)\int W_{s,x,y,z}d\nu 
=(|k|-1)\lan g*\rh_s, \rh_s \ran. 
\edn
Note that
$$
(k - \del_0)+(\check{k} - \del_0)-2(|k|-1)\del_0
=2(|k|-k_0)(p-\del_0),
$$
and that
$$
g*(p-\del_0)=p_{n+1}-\del_0 \ge -\del_0.
$$
Thus, 
\bdnn
\lefteqn{\lan g*(k - \del_{0})*\rh_s, \rh_s \ran 
+\lan g*(\check{k} - \del_{0})*\rh_s, \rh_s \ran 
-2(|k|-1)\lan g*\rh_s, \rh_s \ran } \nn \\
&=& 2(|k|-k_0)\lan g*(p-\del_0)*\rh_s, \rh_s \ran 
\ge 2(|k|-k_0)\cR_s.
\ednn
By this, (\ref{intU}) and (\ref{intW}),  we get
\bdnl{intU>}
\sum_{x,y,z}g(x-y)\int 
\lef( U_{s,x,y,z}-2W_{s,x,y,z}\ri) d\nu 
\ge \lef( \sum_{x,y}g(x-y)\b_{x,y}-2(|k|-k_0) \ri)\cR_s.
\edn
By (\ref{JJ-rr}), (\ref{intV}), (\ref{intU>}) and \Lem{g_n}, 
we obtain (\ref{cA>}) .
\hfill $\Box$

\vvs 
\noindent {\it Proof of part (b)}:
The predictable quadratic variation of the martingale $\cM_\cdot$ 
can be given by:
\bds
\item[1)] \hspace{1cm} 
${\dps 
\lan \cM \ran_t 
= \sum_z \int_0^t ds \int F_{s,z}(\xi)^2\nu(d\xi)}$
\eds
where
$$
F_{s,z}(\xi)=\sum_{x,y} g(x-y) 
( \bar{J}_{s,x,z}(\xi) \bar{J}_{s,y,z}(\xi) - \rho_{s,x} \rho_{s,y} )
$$
Recall that
\bdnn
 \{ \lan \cM \ran_\8 <\8\} & \sub & \{ \mbox{$\cM_t$ converges
as $t \ra \8$}\}\; \; \mbox{a.s.}  \\
\{ \lan \cM \ran_\8 =\8\} & \sub & 
\lef\{ \lim_{t \ra \8}{\cM_t \over \lan \cM \ran_t}=0\ri\}
\; \; \mbox{a.s.} 
\ednn
Thus, to prove (\ref{M=o(intR)}), it is enough to show that 
there is a constant $c \in (0,\8)$ such that
\bds
\item[2)] \hspace{1cm} 
${\dps \lan \cM \ran_t \le c\int_0^t\cR_s ds.}$
\eds
We will do so via two different bounds for $|F_{s,z}(\xi)|$:
\bds
\item[3)] \hspace{1cm} 
$|F_{s,z}(\xi)| \le 2|g|\; \; $ for all $s,z,\xi$,
\item[4)] \hspace{1cm} 
$|F_{s,z}(\xi)| \le c\rh_{s,z}\; \;$ if $\rh_{s,z} \le 1/2$.
\eds 
To get (3), we note that 
$0 \le  \bar{J}_{s,x,z}(\xi) \le 1$ 
and $\sum_x \bar{J}_{s,x,z} = 1$ for each $z$.
Thus, 
\bdmn
|F_{s,z}(\xi)|
& \le &  \lan g*\bar{J}_{s,\cdot,z},\bar{J}_{s,\cdot,z} \ran
+\lan g*\rho_{s},\rho_{s} \ran  \nn \\
& \le & |(g*\bar{J}_{s,\cdot,z})^2|^{1/2} 
|\bar{J}_{s,\cdot,z}^2|^{1/2} + |(g* \rho_s)^2|^{1/2} |\rho_s^2|^{1/2}  \nn \\
& \le & |g||\bar{J}_{s,\cdot,z}^2|+|g|\cR_s \le 2|g|.\nn
\edmn 
To get (4), we assume $\rho_{s,z} \le 1/2$. 
Then, $1+(|\xi|-1)\rho_{s,z} \ge 1/2$ and thus, recalling 
(\ref{Utxyz}) and (\ref{Wtxyz}), 
\bdnn
|F_{s,z}(\xi)| 
&\le & \sum_{x,y}g(x-y)
{|U_{s,x,y,z}(\xi)- W_{s,x,y,z}(\xi)| \over 1+(|\xi|-1)\rho_{s,z}} \\
& \le & 2\sum_{x,y} g(x-y)(|U_{s,x,y,z}(\xi)|+| W_{s,x,y,z}(\xi)|),
\ednn
By (\ref{Wtxyz}) and (\ref{Uexp}), it is clear that the last summation 
is bounded by $c\rh_{s,z}$ for some $c$. \\
(3)--(4) can be used to obtain (2) as follows. 
For each $s$, there is at most one $z$ such that $\rh_{s,z} > 1/2$, 
and $\cR_s >1/4$ if there is such  $z$. 
Thus, 
$$
\sum_z {\bf 1}\{ \rh_{s,z} > 1/2\} <4\cR_s.
$$
By this and (3)--(4), we have 
$$
\sum_z F_{s,z}(\xi)^2 
\le 4|g|^2\sum_z {\bf 1}\{ \rh_{s,z} > 1/2\} 
+c^2\sum_z {\bf 1}\{ \rh_{s,z} \le 1/2\}\rh_{s,z}^2
\le (16|g|^2+c^2)\cR_s.
$$
Plugging this into (1), we are done. 
\hfill $\Box$

\vvs
\small
\noindent{\bf Acknowledgements:}
The authors thank Yuichi Shiozawa for discussions.



\begin{thebibliography}{99}
\bibitem[CH02]{CH02}
Carmona, P., Hu Y.: 
On the partition function of a directed polymer in a  Gaussian random
environment, 
Probab.Theory Related Fields {\bf 124} (2002), no. 3, 431--457. 
\bibitem[CH06]{CH06}
Carmona,P., Hu, Y.:
Strong disorder implies strong localization for directed polymers in a random environment
ALEA Lat. Am. J. Probab. Math. Stat.  {\bf 2} (2006), 217--229.
\bibitem[CSY03]{CSY03}
Comets, F., Shiga, T., Yoshida, N.
Directed polymers in random environment: 
path localization and  strong disorder, 
Bernoulli, {\bf 9} (2003), No. 4, 705--723.
\bibitem[CY05]{CY05}
Comets, F., Yoshida, N.:
Brownian directed polymers in random environment, 
Comm. Math. Phys. {\bf 254} (2005), no. 2, 257--287. 
\bibitem[Gri83]{Gri83}
Griffeath, D.:
The binary contact path process, 
Ann. Probab. {\bf 11} (1983) no. 3, 692--705. 
\bibitem[HWY92]{HWY92}
 Sheng-wu He, Jia-gang Wang, Jia-an Yan:
``Semimartingale theory and stochastic calculus",
Science Press, Beijing New York; 
CRC Press Inc., Boca Raton Ann Arbor London Tokyo, (1992). 
\bibitem[HL81]{HL81}
Holley, R., Liggett, T. M. :
Generalized potlatch and smoothing processes, 
 Z. Wahrsch. Verw. Gebiete {\bf 55} (1981), no. 2, 165--195.
\bibitem[HY09]{HY09}
Hu, Y., Yoshida, N. :
Localization for branching random walks in random environment,
Stochastic Process. Appl.
{\bf 119} (2009), no. 5, 1632--1651.
\bibitem[Lig85]{Lig85}
Liggett, T. M. :
``Interacting Particle Systems",  
Springer Verlag, Berlin-Heidelberg-Tokyo
(1985).
\bibitem[LS81]{LS81}
Liggett, T. M., Spitzer, F. :
Ergodic theorems for coupled random walks and other systems with 
locally interacting components. 
 Z. Wahrsch. Verw. Gebiete {\bf 56} (1981), no. 4, 443--468.
\bibitem[NY09]{NY09} 
Nagahata, Y., Yoshida, N.: Central limit theorem for 
a class of linear systems, 
Electron. J. Probab. {\bf 14} (2009), 960--977.
\bibitem[Sh09]{Sh09} 
Shiozawa, Y.:
Localization for branching Brownian motions in random environment, 
preprint (2009)
\bibitem[Spi81]{Spi81}
Spitzer, F. :
Infinite systems with locally interacting components.  
Ann. Probab. {\bf 9}, (1981), no. 3, 349--364. 
\bibitem[Yo08a]{Yo08a} 
Yoshida, N.:
Phase transitions for the growth rate of linear stochastic evolutions, 
J. Stat. Phys. {\bf 133} (2008), no.6, 1033--1058. 
\bibitem[Yo08b]{Yo08b} 
Yoshida, N.: Localization for linear stochastic evolutions, 
preprint, arXiv:0810.4218, (2008).
\end{thebibliography}
\end{document}